\documentclass[runningheads]{llncs}

\usepackage{graphicx}

\usepackage[short,nocomma]{optidef}
\usepackage{xcolor,soul}
\definecolor{dr}{rgb}{0.75,0.00,0.00}
\definecolor{lr}{rgb}{1.00,0.75,0.75}
\sethlcolor{lr}

\usepackage{amsfonts}
\usepackage[caption=false]{subfig}
\usepackage{booktabs}
\usepackage{threeparttable}
\usepackage{adjustbox}
\usepackage{multirow}
\usepackage{rotating}
\usepackage[lighttt]{lmodern}
\ttfamily
\DeclareFontShape{OT1}{lmtt}{m}{it}
{<->sub*lmtt/m/sl}{}
\usepackage{listings} %https://tex.stackexchange.com/a/234020/191506
\usepackage{opl}
\usepackage{varwidth}
\usepackage[strings]{underscore}
\usepackage[hyphens]{url}
\begin{document}
\title{Optimizing Freight Operations for\\Autonomous Transfer Hub Networks}
\titlerunning{Optimizing Freight Operations for ATHNs}
\author{Kevin Dalmeijer \and Pascal Van Hentenryck}
\authorrunning{K. Dalmeijer and P. Van Hentenryck}
\institute{Georgia Institute of Technology, Atlanta, GA 30332, USA\\
\email{dalmeijer@gatech.edu, pvh@isye.gatech.edu}}

\maketitle              % typeset the header of the contribution

\begin{abstract}
Autonomous trucks are expected to fundamentally transform the freight transportation industry, and the technology is advancing rapidly.
According to some of the major players, the Autonomous Transfer Hub Network (ATHN) business model is the most likely future for the industry.
ATHNs make use of transfer hubs to hand off trailers between human-driven trucks and autonomous trucks.
Autonomous trucks then carry out the transportation between the hubs, while conventional trucks serve the first and last mile.
This paper presents a Constraint Programming (CP) model to schedule the ATHN operations to perform a given set of orders.
The model is used to, for the first time, provide a detailed quantitative study of the benefits of ATHNs by considering a real case study where actual operations are modeled and optimized with high fidelity.
It is found that solving this large-scale optimization problem with CP is computationally feasible, and that ATHNs may bring significant cost savings.
\end{abstract}

\section{Introduction}

Autonomous trucks are expected to fundamentally transform the freight transportation industry, and the enabling technology is progressing rapidly.
Morgan Stanley estimates the potential savings from automation at \$168 billion annually for the US alone \cite{Greene2013-AutonomousFreightVehicles}.
Additionally, autonomous transportation may improve on-road safety, and reduce emissions and traffic congestion \cite{ShortMurray2016-IdentifyingAutonomousVehicle,SlowikSharpe2018-AutomationLongHaul}.

SAE International defines different levels of driving automation, ranging from L0 to L5, corresponding to no-driving automation to full-driving automation \cite{SAEInternational2018-TaxonomyDefinitionsTerms}.
The current focus is on L4 technology (high automation), which aims at delivering automated trucks that can drive without any need for human intervention in specific domains, e.g., on highways.
The automotive industry is actively involved in making L4 vehicles a reality.
Daimler Trucks, one of the leading heavy-duty truck manufacturers in North America, is working with both Torc Robotics and Waymo, and will be testing the latest generation of L4 trucks in the Southwest in early 2021 \cite{Engadget2020-WaymoDaimlerTeam}.
In 2020, truck and engine maker Navistar announced a strategic partnership with technology company TuSimple to develop L4 trucks, to go into production by 2024 \cite{TransportTopics2020-NavistarTusimplePartner}.
Other companies developing self-driving vehicles include Argo AI, Aurora, Cruise, Embark, Ford, Kodiak, Lyft, Motional, Nuro, and Volvo Cars \cite{FleetOwner-TusimpleAutonomousTruck}.

\begin{figure}[!t]
	\centering
	\includegraphics[scale=0.6]{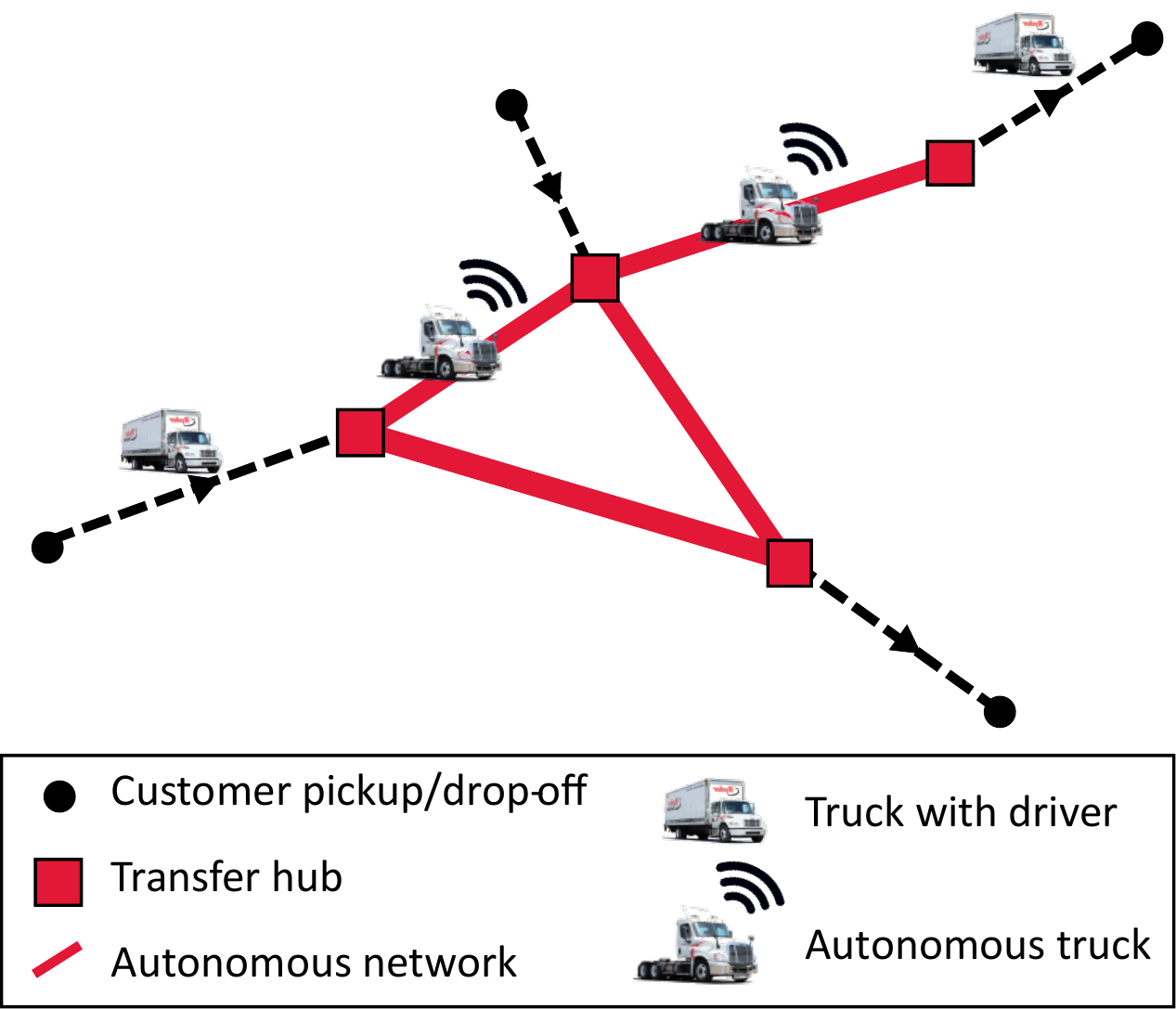}
	\caption{An Example of an Autonomous Transfer Hub Network.}
	\label{fig:autonomous_example}
\end{figure}

A study by Viscelli \cite{Viscelli-Driverless?AutonomousTrucks} describes different scenarios for the adoption of autonomous trucks by the industry.
The most likely scenario, according to some of the major players, is the \emph{transfer hub business model} \cite{Viscelli-Driverless?AutonomousTrucks,RolandBerger2018-ShiftingGearAutomation,ShahandashtEtAl2019-AutonomousVehiclesFreight}.
An Autonomous Transfer Hub Network (ATHN) makes use of autonomous truck ports, or \emph{transfer hubs}, to hand off trailers between human-driven trucks and driverless autonomous trucks.
Autonomous trucks then carry out the transportation between the hubs, while conventional trucks serve the first and last miles.
Figure~\ref{fig:autonomous_example} presents an example of an autonomous network with transfer hubs.
Orders are split into a first-mile leg, an autonomous leg, and a last-mile leg, each of which served by a different vehicle.
A human-driven truck picks up the cargo at the customer location, and drops it off at the nearest transfer hub.
A driverless autonomous truck moves the trailer to the transfer hub closest to the destination, and another human-driven truck performs the last leg.

ATHNs apply automation where it counts: monotonous highway driving is automated, while more complex local driving and customer contact is left to humans.
Especially for long-haul transportation, the benefit of automation is expected to be high.
Global consultancy firm Roland Berger \cite{RolandBerger2018-ShiftingGearAutomation} estimates that operational cost savings may be between 22\% and 40\% in the transfer hub model, based on the cost difference between driverless trucks and conventional trucks.
These estimates are based on single trips and it is not clear that they can be realized in practice: in particular, they do not take into account the empty miles traveled by autonomous trucks to pick up their next orders.

This paper proposes a Constraint Programming (CP) model to schedule the ATHN operations for a given set of orders.
The resulting schedule details the autonomous operations and the first/last-mile operations at each of the hubs, and specifies the movements of every load, vehicle, and driver.
\emph{The CP model is then used to provide, for the first time, a detailed quantitative study of the benefits of ATHNs by considering a real case study where actual operations are modeled and optimized with high fidelity.}
It examines whether the savings predicted by \cite{RolandBerger2018-ShiftingGearAutomation} materialize when the network effects, e.g., empty miles for relocation, are taken into account.
It is found that it is computationally feasible to solve this large-scale optimization problem with more than 100,000 variables and 100,000 constraints, and that the benefits of ATHNs may indeed be realized in practice.

The remainder of this paper is organized as follows.
Section~\ref{sec:problem} defines the problem of scheduling freight operations on an ATHN, and Section~\ref{sec:formulation} formulates a CP model to solve this problem.
The case study is presented in Section~\ref{sec:casestudy}.
The final section of the paper summarizes the findings and provides the conclusions.

\section{Problem Statement}
\label{sec:problem}

This section defines the problem of scheduling freight operations on an ATHN to perform a given set orders with a given set of vehicles, with the objective to minimize the cost of driving empty.
The problem is defined on a directed graph $G=(V,A)$, with vertices $V$ and arcs $A$.
The vertices represent locations, and are partitioned into hub locations $V_H$ and customer locations $V_C$.
Arcs between the transfer hubs correspond to autonomous transportation, and the other arcs represent human-driven legs.
Every arc $(i,j) \in A$ is associated with a non-negative travel time $\tau_{ij}$ and a cost $c_{ij}$.
For convenience, define $\tau_{ii} = 0$ and $c_{ii} = 0$ for all $i \in V$.

The set of customer orders is given by $R$.
Order $r\in R$ is supposed to be picked up at time $p(r)$ at the origin $o(r) \in V_C$, and to be transported through the ATHN to the destination $d(r) \in V_C$.
Based on order $r\in R$, three tasks are defined: the first-mile task $t_r^f$, the autonomous task $t_r^a$, and the last-mile task $t_r^l$.
The first-mile task consists of loading the trailer at the customer location, moving the freight to the closest transfer hub, and unloading the trailer.
Similarly, the autonomous task and the last-mile task consist of loading, driving (between the hubs and to the destination, respectively), and unloading.

Let $T$ be the set of all tasks generated by the orders.
Every task $t\in T$ corresponds to a single leg, and is defined by an origin $o(t) \in V$, a destination $d(t) \in V$, and a pickup time $p(t)$.
The duration of a task equals $\tau_{o(t),d(t)} + 2S$, where $S \ge 0$ is the fixed time for loading or unloading a trailer.
The pickup time $p(t_r^f)$ of the first-mile task is equal to the order pickup time $p(r)$, while subsequent pickup times are based on the time the freight is supposed to be available.
That is, $p(t_r^a) = p(t_r^f) + \tau_{o(t_r^f),d(t_r^f)} + 2S$, and $p(t_r^l) = p(t_r^a) + \tau_{o(t_r^a),d(t_r^a)} + 2S$.

To create a feasible schedule, every task must be given a starting time, and be assigned to one of the available trucks.
It is assumed that an appointment flexibility of $\Delta \ge 0$ minutes is permitted, which means that task $t\in T$ may start anywhere in the interval $[p(t)-\Delta, p(t)+\Delta]$.
The set of trucks $K$ is partitioned into autonomous trucks $K_A$, and regular trucks $K_h$ at every hub $h \in V_H$.
Tasks can only be assigned to the corresponding set of trucks, and tasks performed by the same vehicle must not overlap in time.
If $t\in T$ and $t'\in T$ are subsequent tasks for a single truck, and $d(t) \neq o(t')$, then an empty relocation is necessary, which takes $\tau_{d(t) o(t')}$ time units and has cost $c_{d(t) o(t')}$.
The objective is to assign the tasks such that the total relocation cost is minimized.

Note that the problem described above can be decomposed and solved independently for the autonomous network and for the operations at each of the hubs.
This is possible because different trucks are used for each part of the ATHN, and because each task is given an independent pickup time, based on the expected time the freight is available.
One potential problem is that using the appointment flexibility in one part of the network may lead to an infeasibility in another part of the network, but the case study shows that this is not an issue in practice: The first and last-mile schedules are not very constrained, and shifting the schedule to accommodate flexibility in the autonomous network is straightforward.
Alternatively, one may first optimize the autonomous network, and update the first and last-mile pickup times accordingly.

\section{Mathematical Formulation}
\label{sec:formulation}

This section presents a CP model to schedule orders on the ATHN.
Without loss of generality, the set of tasks and the set of trucks represent a single part of the network that can be optimized independently.
That is, either the autonomous operations, or the first/last-mile operations at one of the hubs are considered.

\newsavebox{\modelbox}
\begin{lrbox}{\modelbox}
\begin{varwidth}{1.15\textwidth}
\begin{lstlisting}
range Trucks = ...;
range Tasks = ...;
range Sites = ...;
range Horizon = ...;
range Types = Sites union { shipType };  
int or[Tasks] = ...;                     
int de[Tasks] = ...;                     
int pickupTime[Tasks] = ...;
int loadTime = ...;
int flexibility = ...;
int travelTime[Types,Types] = ...;
int travelCost[Types,Types] = ...;       

dvar interval task[t in Tasks] in Horizon
              size travelTime[or[t],de[t]] + 2*loadTime;
dvar interval ttask[k in Trucks,t in Tasks] optional in Horizon
              size travelTime[or[t],de[t]] + 2*loadTime;
dvar interval load[Trucks,Tasks] optional in Horizon size loadTime;
dvar interval ship[k in Trucks,t in Tasks] optional in Horizon
              size travelTime[ort],de[t]];
dvar interval unload[Trucks,Tasks] optional in Horizon size loadTime;
dvar sequence truckSeq[k in Trucks]
  in append(all(t in Tasks)load[k,t],all(t in Tasks)ship[k,t],all(t in Tasks)unload[k,t])
  types append(all(t in Tasks)or[t],all(t in Tasks)shipType,all(t in Tasks)de[t]);
dvar int emptyMilesCost[Trucks,Tasks];
dvar int truckEmptyMilesCost[Trucks];

minimize sum(k in Trucks) truckEmptyMilesCost[k];

constraints {

   forall(t in Tasks) 
      startOf(task[t]) >= pickupTime[t] - flexibility;
      startOf(task[t]) <= pickupTime[t] + flexibility;
	
   forall(k in Trucks,t in Tasks)
      span(ttask[k,t],[load[k,t],ship[k,t],unload[k,t]]);
      startOf(ship[k,t]) == endOf(load[k,t])
      startOf(unload[k,t]) == endOf(ship[k,t])	 
	
   forall(k in Trucks)
      alternative(task[t],all(k in Trucks) ttask[k,t])	
	
   forall(k in Trucks,t in Tasks)
      emptyMilesCost[k,t] = travelCost[destination[t],typeOfNext(truckSeq[k],ttask[k,t],destination[t],destination[t])];
	
   forall(k in Trucks)
      truckEmptyMilesCost[k] = sum(t in Tasks) emptyMilesCost[k,t];
	
   forall(k in Trucks)
      noOverlap(truckSeq,travelTime);

}
\end{lstlisting}
\end{varwidth}
\end{lrbox}

\begin{figure}[!t]
\makebox[\textwidth][c]{%
\fbox{\begin{minipage}{1.13\textwidth}
	\usebox{\modelbox}
\end{minipage}}
}
\caption{Formulation for Scheduling Freight Operations on an ATHN.}
\label{fig:formulation}
\end{figure}

The model is depicted in Figure \ref{fig:formulation} using OPL
syntax \cite{VanHentenryck1999-OplOptimizationProgramming}. The data of the model is given in lines 1--12. It consists of
a number of ranges (line 1--5), information about the tasks (lines
6--8) that include their origins, destinations, and pickup times, the
time to load/unload a trailer (line 9), the flexibility around the
pickup times (line 10), and the matrices of travel times and travel
costs. These matrices are defined between the sites but also
include a dummy location {\tt shipType} for reasons that will become
clear shortly.

The main decision variables are the interval variables {\tt task[t]}
that specify the start and end times of task {\tt t} when processed
by the autonomous network, and the optional interval variables {\tt
	ttask[k,t]} that are present if task {\tt t} is transported by
truck {\tt k}. These optional variables consist of three subtasks that
are captured by the interval variables {\tt load[k,t]} for loading,
{\tt ship[k,t]} for transportation, and {\tt unload[k,t]} for
unloading. The other key decision variables are the sequence variables
{\tt truckSeq[k]} associated with every truck: these variables
represent the sequence of tasks performed by every truck. They
contain the loading, shipping, and unloading interval variables
associated with the trucks, and their types. The type of a loading
interval variable is the origin of the task, the type of an unloading
interval variable is the destination of the task, and the type of the
shipping interval variable is the specific type {\tt shipType} that is
used to represent the fact that there is no transition cost and
transition time between the load and shipping subtasks, and the
shipping and destination subtasks. The model also contains two
auxiliary decision variables to capture the empty mile cost between a
task and its successor, and the empty mile cost of the truck
sequence.

The objective function (line 28) minimizes the total costs of empty
miles. The constraints in lines 32--34 specify the potential start
times of the tasks, and are defined in terms of the pickup times and
the flexibility parameter. The {\sc span} constraints (line 37) link
the task variables and their subtasks, while the constraints in lines
38--39 link the subtasks together. The {\sc alternative} constraints
on line 42 specify that each task is processed by a single truck. The
empty mile costs between a task and its subsequent task (if it
exists) is computed by the constraints in line 45: they use the {\sc
	typeOfNext} expression on the sequence variables. The total empty
mile cost for a truck is computed in line 48. The {\sc noOverlap}
constraints in line 51 impose the disjunctive constraints between the
tasks and the transition times.

\section{Case Study}
\label{sec:casestudy}

To quantify the impact of autonomous trucking on a real transportation network, a case study is presented for the dedicated transportation business of Ryder System, Inc., commonly referred to as \emph{Ryder}.
Ryder is one of the largest transportation and logistics companies in North America, and provides fleet management, supply chain, and dedicated transportation services to over 50,000 customers.
Its dedicated business, \emph{Ryder Dedicated Transportation Solutions}, offers supply-chain solutions in which Ryder provides both drivers and trucks, and handles all other aspects of managing the fleet.
Ryder's order data is used to design an ATHN, and to create a detailed plan for how it would operate.
This allows for a realistic evaluation of the benefits of autonomous trucking.

\subsection{Data Description}
\label{sec:inputdata}

Ryder prepared a representative dataset for its dedicated transportation business in the Southeast of the US, reducing the scope to orders that were strong candidates for automation.
The dataset consists of trips that start in the first week of October 2019, and stay completely within the following states: AL, FL, GA, MS, NC, SC, and TN.
It contains 11,264 rows, which corresponds to 2,090 orders, formatted as in Table~\ref{tab:order_data}.
Every order has a unique \emph{OrderNumber}, and every row corresponds to a stop for a particular order.
Stops have a unique identifier \emph{StopNumber}, and the \emph{Stop} column indicates the sequence within the order.
The columns \emph{StopArrivalDate} and \emph{StopDepartureDate} indicate the scheduled arrival and departure times, and \emph{City} and \emph{ZipCode} identify the location of the stop.
The \emph{Status} column gives a code for the status of the vehicle on arrival, and the \emph{Event} column indicates what happens at the stop.

\begin{adjustbox}{center,float={table}[!t]}
	\centering
	\scriptsize
	\begin{threeparttable}
		\caption{An Example of the Order Data.}
		\label{tab:order_data}%
		\begin{tabular}{rrrrrrrrrr}
			\multicolumn{1}{l}{StopNum} & \multicolumn{1}{l}{OrderNum} & \multicolumn{1}{l}{StopArrivalDate} & \multicolumn{1}{l}{StopDepartureDate} & \multicolumn{1}{l}{Stop} & \multicolumn{1}{l}{City} & ZipCode & \multicolumn{1}{l}{Status} & Event &  \\
			\toprule
			68315760 & 7366366 & 2-10-2019 09:01 & 2-10-2019 09:02 & 1     & Atlanta & 30303 & LD    & HPL \\
			68315761 & 7366366 & 2-10-2019 16:29 & 2-10-2019 18:33 & 2     & Tennessee & 37774 & LD    & LUL \\
			68315762 & 7366366 & 3-10-2019 11:00 & 3-10-2019 11:30 & 3     & Atlanta & 30303 & MT    & DMT \\
			% 31301271 & 8846813 & 2-10-2019 08:30 & 2-10-2019 08:45 & 1     & Savannah & \fix{00000} & NaN   & LLD \\
			% 31301272 & 8846813 & 2-10-2019 14:45 & 2-10-2019 16:00 & 2     & Jacksonville & \fix{00000} & LD    & LUL \\
			\dots & \dots & \dots & \dots & \dots & \dots & \dots & \dots & \dots\\
			46798427 & 5207334 & 7-10-2019 02:35 & 7-10-2019 02:50 & 1     & Alpharetta & 30009 & NaN   & LLD \\
			46798428 & 5207334 & 7-10-2019 08:10 & 7-10-2019 08:49 & 2     & Macon & 31201 & LD    & LUL \\
			46798429 & 5207334 & 7-10-2019 15:16 & 7-10-2019 15:45 & 3     & Alpharetta & 30009 & LD    & LUL \\
		\end{tabular}%
	\end{threeparttable}
\end{adjustbox}%

The example data in Table~\ref{tab:order_data} displays two orders.
The first order is a trip from Atlanta to Tennessee and back.
Based on the status code, the truck arrived in Tennessee loaded (LD) and returned to Atlanta empty (MT).
The event codes show that a preloaded trailer was hooked in Atlanta (HPL), followed by a live unload (LUL) in Tennessee, after which the truck dropped the empty trailer (DMT) in Atlanta.
The exact codes are not important for the purpose of this paper.
What is important, is the ability to derive the parts of the trip when the truck is moving freight, and the parts when the truck is driving empty.
If a vehicle returns to the starting location after only making deliveries, it is assumed that the return trip is empty, and the data is corrected if needed.

Road system data was obtained from OpenStreetMap \cite{OpenStreetMap2020}, and route distance and mileage were calculated with the GraphHopper library.
The provided orders are long-haul trips, with an average trip length of 431 miles.
Most of this distance is driven on highways: 65\% of the distance is driven on interstates and US highways, and this number goes up further to 87\% if state highways are included.
Significant highway usage is typical for long-haul transportation, and indicates that a significant part of each trip can potentially be automated.

\subsection{ATHN Design}

\begin{figure}[!t]
	\centering
	\subfloat[Small Network (17 hubs).\label{fig:designsmall}]{%
		\centering
		\includegraphics[width=0.45\textwidth, trim=28cm 3cm 29cm 15cm, clip]{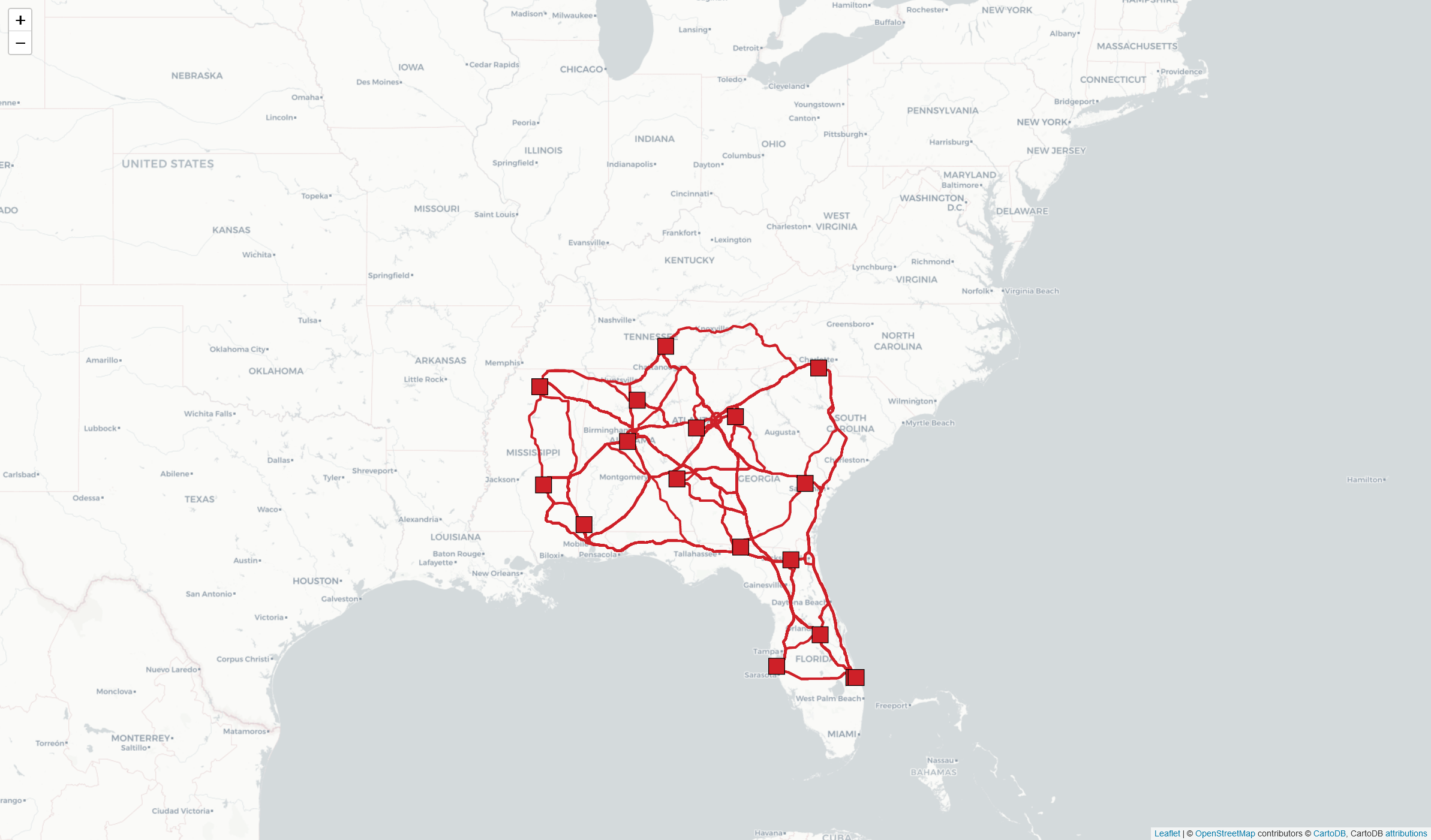}
	}
	\hfill
	\subfloat[Large Network (30 hubs).\label{fig:designlarge}]{%
		\centering
		\includegraphics[width=0.45\textwidth, trim=28cm 3cm 29cm 15cm, clip]{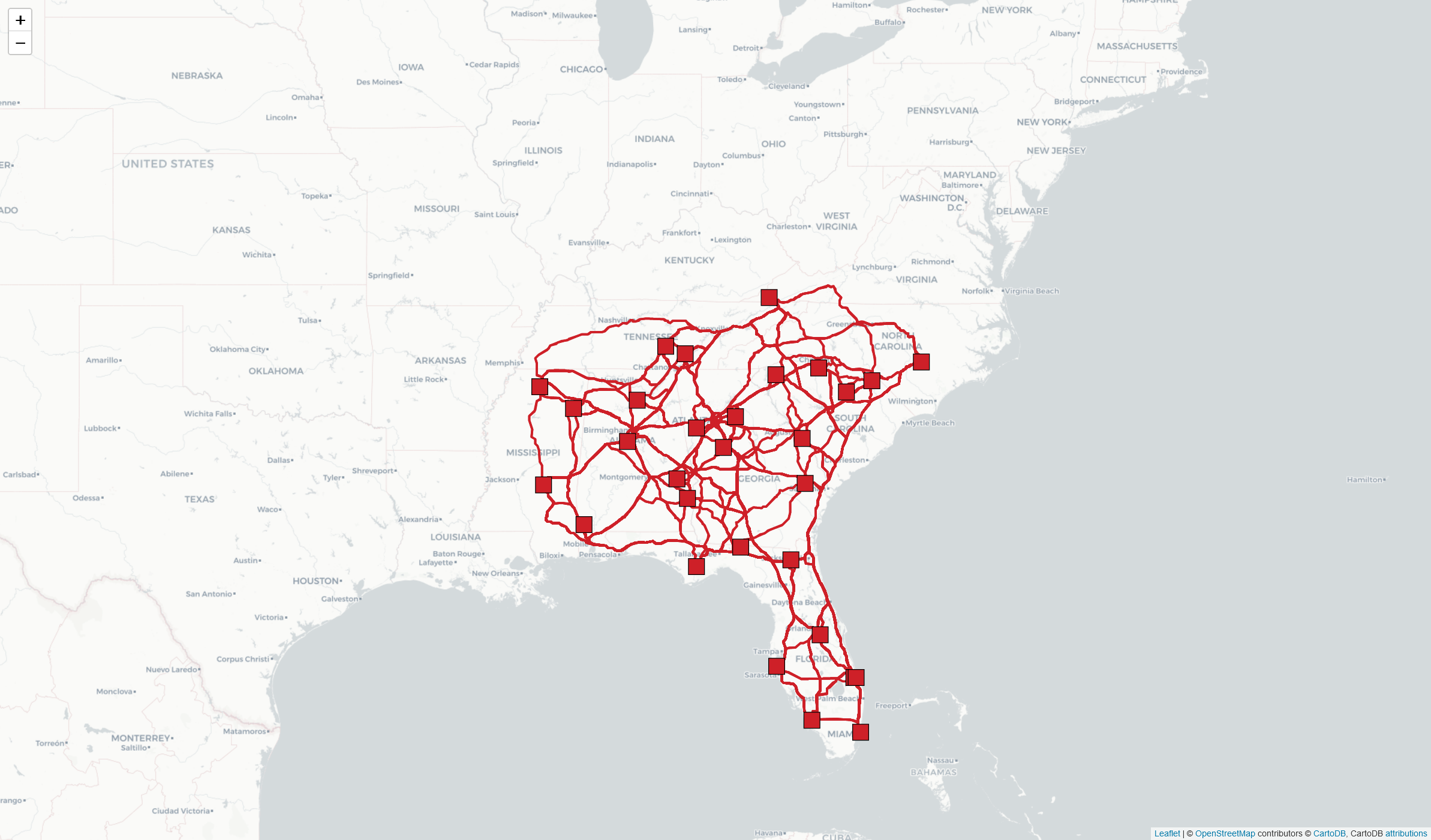}
	}
	\caption{ATHN Network Designs for the Southeast.}
	\label{fig:designs}
\end{figure}

The design of ATHNs needs to decide the locations of the transfer hubs, which are the gateways to the autonomous parts of the network.
A natural choice is to locate the hubs in areas where many trucks currently enter or exit the highway system.
Historical order data is used to identify these common highway access points.
For a given order, the truck is routed through the existing road network, and the highway segments, and their access points, can easily be identified.
The transfer hubs are then placed in areas with many access points.
The case study considers two different sets of hubs: a \emph{small network} with 17 transfer hubs in the areas where Ryder trucks most frequently access the highway system, and a \emph{large network} that includes 13 additional hubs in locations with fewer highway access points.
The small and the large network are visualized in Figure~\ref{fig:designsmall} and Figure~\ref{fig:designlarge}, respectively.
The exact hub locations are masked, but the figures are accurate within a 50 mile range.
The large network extends further northeast into North-Carolina, and further south into Florida.
It also makes the network more dense in the center of the region.

\subsection{Order Selection}

The case study focuses on scheduling the 494 most \emph{challenging orders}.
These orders consist of a single delivery, followed by an empty return trip.
Because they require empty travel for 50\% of the trip, these challenging orders are an excellent target for cost savings.
They also make up 24\% of the dataset, and account for 53\% of the empty mileage.

For a given ATHN and given arc costs $c_{ij}$ for all $(i,j) \in A$, it is first determined which of these orders may benefit from the autonomous network, and which are better served by human-driven trucks.
For example, a direct trip with a human-driven trucks may be preferred for trips that are short, or for trips between locations that are far from any hub.
Two options are compared for every order: The first option is to serve the order with a conventional truck, which necessitates an empty return leg.
The second option uses the autonomous network, which amounts to driving to the nearest hub with a human-driven truck (first mile), shipping the load over the autonomous network to the hub closest to the destination, and then serving the last mile with another human-driven truck.
No empty returns are added in this case, as every truck is immediately available for the next task.
The two options are compared in terms of arc costs, and the cheapest one is selected.
Different costs are used for autonomous and non-autonomous arcs, as will be explained in the next section.
The CP model only considers the orders that may benefit from the ATHN, while the other orders are scheduled separately.

\subsection{Parameters and Settings}

The following scenario is defined as the \emph{base case} for the upcoming experiments.
The base case uses the small network, presented by Figure~\ref{fig:designsmall}.
For conventional trucks, the cost $c_{ij}$ for driving arc $(i,j)\in A$ is equal to the road distance.
For autonomous trucks, this cost is reduced by the percentage $\alpha$.
The value of $\alpha$ is a parameter, and the case study uses values ranging from $\alpha=25\%$ to $\alpha=40\%$, with $\alpha = 25\%$ for the base case.
This results in a conservative estimate of the benefits of autonomous trucking, which is predicted to be 29\% to 45\% cheaper per mile \cite{EngholmEtAl2020-CostAnalysisDriverless}.

The time for loading or unloading is estimated at $S = 30$ minutes, and the appointment time flexibility is set to one hour ($\Delta=60$).
The number of autonomous trucks is set to $\lvert K_A \rvert=50$.
The CP model is used to schedule the orders for each independent part of the network: the autonomous operations, and the first/last-mile operations at each of the hubs.
Each model is solved with the CPLEX CP Optimizer version 12.8 \cite{LaborieEtAl2018-IbmIlogCp}.

\subsection{Base Case Results}

\begin{figure}[p]
	\centering
	\includegraphics[trim=20 50 10 40,clip,width=\linewidth]{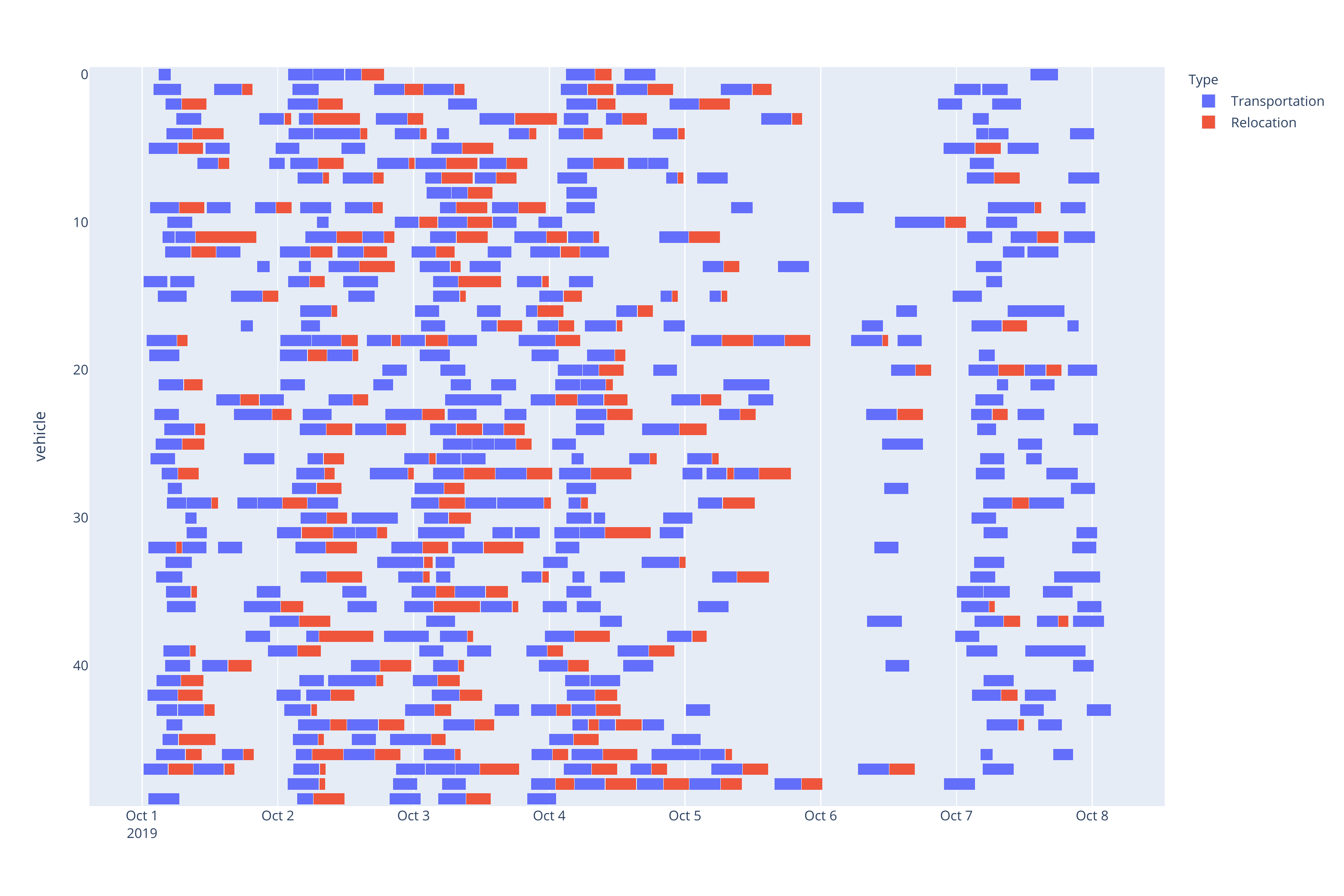}
    \caption{Autonomous Truck Schedule for the Base Case.}
	\label{fig:basecase_schedule}
\end{figure}

Out of the 494 challenging orders, 437 (88\%) are found to potentially benefit from the autonomous network, and the CP model is used to schedule these orders.
Scheduling the autonomous part of the ATHN is the most challenging: the model has close to 110,000 decision variables and more than 110,000 constraints.
The model is given an hour of CPU time and returns the best found solution, which is visualized in Figure~\ref{fig:basecase_schedule}.
This figure shows both the autonomous tasks and the and the relocation tasks for the first week of October.

\begin{figure}[p]
	\centering
	\includegraphics[width=0.7\linewidth]{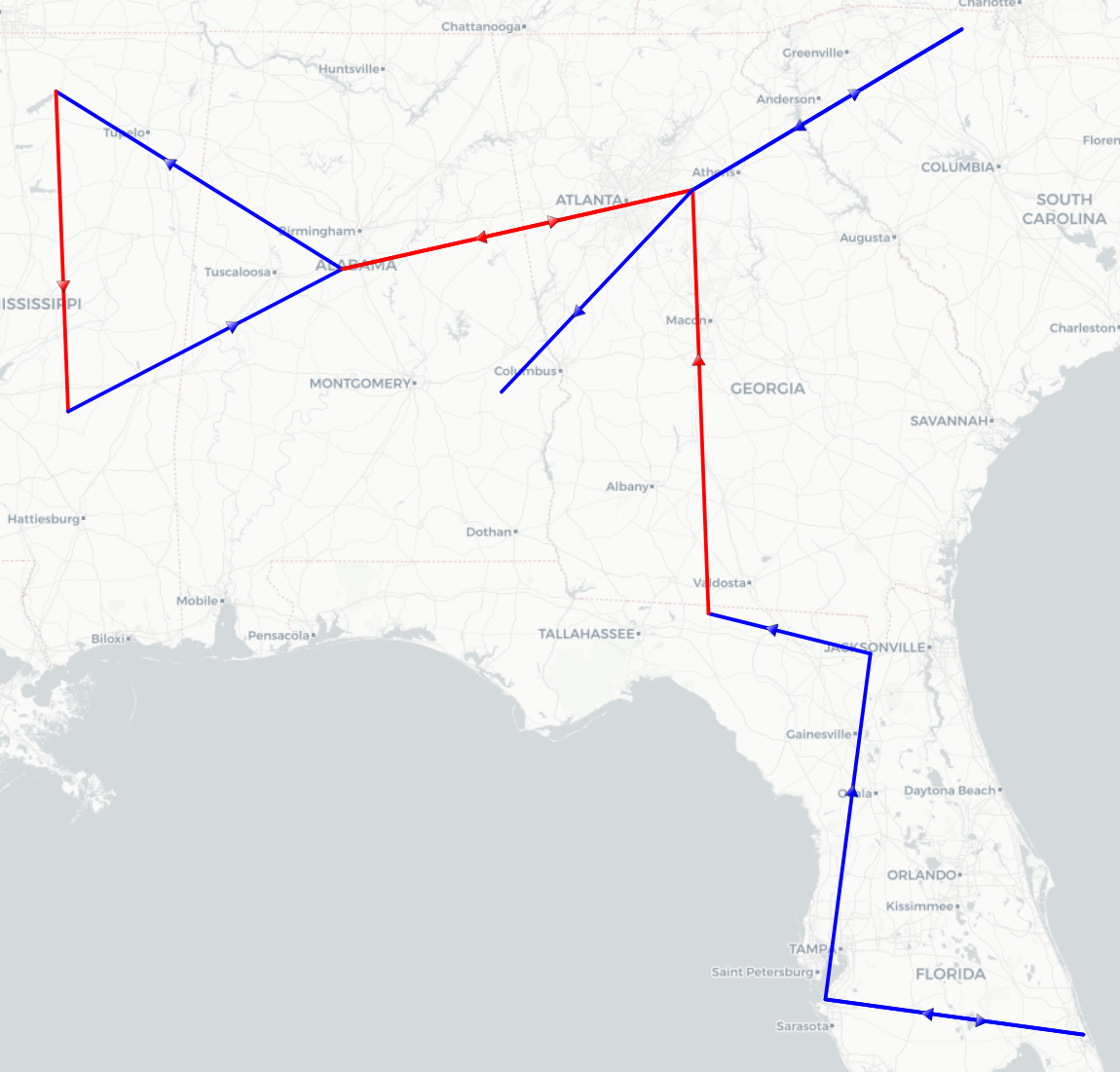}
	\caption{Single Autonomous Truck Route in the Base Case (blue is loaded, red is empty).}
	\label{fig:singleroute}
\end{figure}

Figure~\ref{fig:basecase_schedule} shows that the transportation tasks are close together, and only a relatively small amount of relocation is necessary.
It is interesting to see that only a small number of autonomous trucks is driving during the weekend.
This is because the appointment times are still based on the current agreements with the customers, and having drivers work during the weekends is typically avoided.
For autonomous trucks, this would not be a problem, which again underlines that making full use of autonomous transportation requires adapting the business model and current practices.
In Section~\ref{sec:flexibility}, the importance of time flexibility is considered in more detail.

The routes that are driven by the autonomous trucks consist of serving autonomous legs of different orders, with relocations in between.
Figure~\ref{fig:singleroute} shows a representative single truck route from the schedule.
Blue arrows indicate that freight is being moved, and red arrows indicate that the truck is driving empty.
The truck starts at the west coast of Florida, where it picks up a load that has been delivered to the hub by a first/last-mile truck driver.
The freight is then transported to the east coast, where it is unloaded so that a regular truck with driver can complete the last-mile.
The autonomous truck immediately starts serving the autonomous leg of the next order, returning to the Tampa area.
After that, legs are served that are going north.
The first time a relocation is needed, is when the truck makes a delivery near Valdosta, close to the Georgia and Florida border.
No freight is immediately available, and the vehicle drives empty to the next location to pick up a load there.
The routes are clearly complex, which emphasizes the power of optimization: It is unlikely that this solution can be found by manual planners, but it is possible with optimization techniques.

Compared to the autonomous trucks, the total distance driven by regular trucks is relatively short.
The optimization model was used to schedule the first and last-mile operations at selected hubs, and it was found that the amount of work is often insufficient to keep drivers occupied throughout the week.
To prevent driver idle time, it may be beneficial to outsource these legs, or to consolidate them with other operations in the area.
In terms of mileage, the percentage of empty miles at the hubs is typically under 25\%, which is used as an estimate for the first/last-mile efficiency in the remainder.

Table~\ref{tab:basecase_costs} quantifies the impact of autonomous trucking on the operating costs for the 437 selected orders.
The \emph{Mileage} column indicates the total miles driven for both the current network and for the ATHN.
The numbers are separated based on whether the distance was driven while loaded or empty, and percentages are shown in the \emph{\% of total} column.
The \emph{Cost without autonomous trucks} converts the miles into dollars, using \$2 per mile as an approximation for the cost of human-driven trucks.
Recall that driving autonomously is assumed to be $\alpha=25\%$ cheaper, which is reflected by the \emph{Cost adjustment} column.
The \emph{Cost} column presents the cost when autonomous trucks are available, and is obtained by multiplying the cost without autonomous trucks by the cost adjustment factor.

\begin{adjustbox}{center,float={table}[!t]}
	\centering
	\footnotesize
	\begin{threeparttable}
		\caption{Cost Table for the Base Case (437 orders).}
		\label{tab:basecase_costs}%
		\begin{tabular}{ccrrrrrrr}
			\toprule
			&		  &		  &		    &			  & Cost without  &           & \\
			&       &       & \quad Mileage & \quad \% of total & \quad auton. trucks & \quad Cost adj. & \quad\quad\quad Cost \\
			\midrule
			\multirow{3}[0]{*}{Current network} &       & \multicolumn{1}{l}{Loaded} & 96,669 & 50\%  & \$ 193,338 & 1.00  & \$ 193,338 \\
			&       & \multicolumn{1}{l}{Empty} & 96,698 & 50\%  & \$ 193,396 & 1.00  & \$ 193,396 \\
			&       & \multicolumn{1}{l}{Total} & 193,367 & 100\% & \$ 386,734 & 1.00  & \$ 386,734 \\
			\midrule
			\multirow{9}[0]{*}{\parbox{2.5 cm}{Autonomous transfer\\hub network}} & \multirow{3}[0]{*}{Autonomous} & \multicolumn{1}{l}{Loaded} & 91,618 & 67\%  & \$ 183,235 & 0.75  & \$ 137,426 \\
			&       & \multicolumn{1}{l}{Empty} & 44,217 & 33\%  & \$ 88,433 & 0.75  & \$ 66,325 \\
			&       & \multicolumn{1}{l}{Total} & 135,834 & 100\% & \$ 271,668 & 0.75  & \$ 203,751 \\
			\cmidrule{2-8}
			& \multirow{3}[0]{*}{First/last mile} & \multicolumn{1}{l}{Loaded} & 29,286 & 75\%  & \$ 58,573 & 1.00  & \$ 58,573 \\
			&       & \multicolumn{1}{l}{Empty\tnote{*}} & 9,762 & 25\%  & \$ 19,524 & 1.00  & \$ 19,524 \\
			&       & \multicolumn{1}{l}{Total} & 39,049 & 100\% & \$ 78,097 & 1.00  & \$ 78,097 \\
			\cmidrule{2-8}
			& Total &       & 174,883 &       & \$ 349,766 &       & \$ 281,848 \\
			\cmidrule{2-8}
			& Savings &       & 18,484 &       & \$ 36,969 &       & \$ 104,886 \\
			& Savings (\%) &       & 10\% &       & 10\% &       & 27\% \\
			\bottomrule
		\end{tabular}%
		\begin{tablenotes}
			\item[*] estimated
		\end{tablenotes}
	\end{threeparttable}
\end{adjustbox}%

Compared to the current network, the ATHN allows for significant savings for the selected orders: Table~\ref{tab:basecase_costs} shows that the total cost goes down by 27\%.
At a cost of \$2 per mile, this corresponds to \$104,886 per week, or \$5.5M per year.
The \emph{Mileage} column shows that almost 80\% of the mileage in the ATHN can be automated, which partly explains the large savings.
What is very interesting to observe is that the total mileage for the ATHN is actually \emph{less} than the total mileage for the direct trips in the current network.
In the transfer hub network, there is no need to return back empty after a delivery, and there is no need to limit working hours or to return to a domicile at the end of the day.
As a result, only 33\% of the automated distance is driven empty, compared to 50\% for the current system.
This means that even if autonomous trucks would be as expensive as trucks with drivers, costs would still go down by 10\% due to the additional flexibility that automation brings.

\subsection{Impact of the Size of the Network}

A larger autonomous network results in shorter first/last-mile trips, and may have a larger area of coverage.
To evaluate the impact of the size of the network, the calculations for the base case are repeated using the large network (Figure~\ref{fig:designlarge}) with 30 hubs, instead of the small network (Figure~\ref{fig:designsmall}) with 17 hubs.
For the large network, 468 of the 494 orders (95\%) may benefit from the autonomous network, compared to only 88\% for the base case.
This immediately implies that there is more potential for savings.
It also means a higher utilization of the autonomous trucks, as more legs are served by the same 50 vehicles.

Table~\ref{tab:large_25_costs} shows that the relative cost savings for the large network (29\%) are similar to those for the small network (27\%).
This means that the average benefit of automation is similar for both designs, \emph{for the orders that are automated}.
However, the large network allows more trips to benefit from automation, which is why the cost savings of \$ 116,582 are 11\% higher than the savings for the small network (\$ 104,886).
The average benefit of automation is similar for the two designs due to two effects that cancel out.
First, the same autonomous trucks have to serve more orders on the large network.
This increases the utilization of the vehicles, but also increases the percentage of empty miles from 33\% to 35\%.
The reason for this increase is that there is less time available to wait around at a hub for the next order, as the trucks are needed to perform other orders in the meantime.
On the other hand, the first and last-mile trips are shorter due to the additional hubs, which saves costs.

\begin{adjustbox}{center,float={table}[!t]}
	\centering
	\footnotesize
	\begin{threeparttable}
		\caption{Cost Table for the Large Network (468 orders).}
		\label{tab:large_25_costs}%
		\begin{tabular}{ccrrrrrr}
	    \toprule
		&		  &		  &		    &			  & Cost without  &           & \\
		&       &       & \quad Mileage & \quad \% of total & \quad auton. trucks & \quad Cost adj. & \quad\quad\quad Cost \\
		\midrule
		\multirow{3}[0]{*}{Current network} &       & \multicolumn{1}{l}{Loaded} & 101,213 & 50\%  & \$ 202,425 & 1.00  & \$ 202,425 \\
		&       & \multicolumn{1}{l}{Empty} & 96,698 & 50\%  & \$ 193,396 & 1.00  & \$ 193,396 \\
		&       & \multicolumn{1}{l}{Total} & 202,476 & 100\% & \$ 404,953 & 1.00  & \$ 404,953 \\
		\midrule
		\multirow{8}[0]{*}{\parbox{2.5 cm}{Autonomous transfer\\hub network}} & \multirow{3}[0]{*}{Autonomous} & \multicolumn{1}{l}{Loaded} & 97,326 & 65\%  & \$ 194,653 & 0.75  & \$ 145,990 \\
		&       & \multicolumn{1}{l}{Empty} & 53,247 & 35\%  & \$ 106,493 & 0.75  & \$ 79,870 \\
		&       & \multicolumn{1}{l}{Total} & 150,573 & 100\% & \$ 301,146 & 0.75  & \$ 225,860 \\
		\cmidrule{2-8}
		& \multirow{3}[0]{*}{First/last mile} & \multicolumn{1}{l}{Loaded} & 23,442 & 75\%  & \$ 46,883 & 1.00  & \$ 46,883 \\
		&       & \multicolumn{1}{l}{Empty \tnote{*}} & 7,814 & 25\%  & \$ 15,628 & 1.00  & \$ 15,628 \\
		&       & \multicolumn{1}{l}{Total} & 31,256 & 100\% & \$ 62,511 & 1.00  & \$ 62,511 \\
		\cmidrule{2-8}
		& Total &       & 181,829 &       & \$ 363,657 &       & \$ 288,371 \\
		\cmidrule{2-8}
		& Savings &       & 20,648 &       & \$ 41,296 &       & \$ 116,582 \\
		& Savings (\%) &       & 10\%  &       & 10\%  &       & 29\% \\
		\bottomrule
		\end{tabular}%
		\begin{tablenotes}
			\item[*] estimated
		\end{tablenotes}
	\end{threeparttable}
\end{adjustbox}%

\subsection{Impact of the Cost of Autonomous Trucking}

For the base case, it was assumed that autonomous trucks are $\alpha= 25\%$ cheaper per mile than trucks with a driver.
However, this number is yet far from certain, and higher cost reductions have also been reported in the literature.
To investigate the impact of the cost of autonomous trucking, Table~\ref{tab:overview_costs} presents results for $\alpha$ ranging from $25\%$ to $40\%$, for both the small and the large network.
The column \emph{Autom.
	orders} gives the number of orders that may benefit from automation, and are considered in the ATHN.
The relative cost savings (\emph{Rel.
	savings}) state the cost reduction compared to serving these orders with conventional trucks.
The \emph{Cost savings} column gives the absolute cost savings in dollars.
The final column compares the absolute savings to the savings obtained for the baseline (small network, $\alpha=25\%$).

Table~\ref{tab:overview_costs} shows that, as autonomous trucking gets cheaper, and as more hubs are added to the network, the savings compared to the current system go up.
Additionally, more orders start using the ATHN, which increases the absolute savings further.
Even though the autonomous trucks only perform the transportation between the hubs, the relative cost savings for the complete system often exceed the mileage cost reduction for autonomous trucks ($\alpha$).
This again shows that the benefit of autonomous trucks is not only the lower cost per mile, but also the additional flexibility.
Compared to the base case, cheaper autonomous transportation results in significantly larger savings.
Similar as in the previous section, increasing the size of the network does not strongly impact the average cost benefit per order, but does increase the total amount of orders that can be automated, which leads to more profits.
In the best case (large network, $\alpha=40\%$), the total benefit of the ATHN is \$ 161,762 per week for the challenging orders, which corresponds to \$8.4M savings per year.

\begin{table}[!t]
	\centering
	\footnotesize
	\caption{Overview of Cost Savings under Different Assumptions.}
	\label{tab:overview_costs}%
	\begin{tabular}{crrrrr}
		\toprule
		&       &       &       &       & Additional savings \\
		\multicolumn{1}{l}{Network} & $\alpha$ & \quad Autom. orders & \quad Rel. savings & \quad Cost savings & \quad comp. to base case \\
		\midrule
		\multirow{4}[0]{*}{Small} & 25\%  & 437   & 27\%  & \$ 104,886 & +0\% \\
		& 30\%  & 439   & 32\%  & \$ 122,396 & +17\% \\
		& 35\%  & 443   & 35\%  & \$ 135,572 & +29\% \\
		& 40\%  & 443   & 38\%  & \$ 148,984 & +42\% \\
		\midrule
		\multirow{4}[0]{*}{Large} & 25\%  & 468   & 29\%  & \$ 116,582 & +11\% \\
		& 30\%  & 469   & 33\%  & \$ 131,798 & +26\% \\
		& 35\%  & 472   & 37\%  & \$ 149,586 & +43\% \\
		& 40\%  & 472   & 40\%  & \$ 161,762 & +54\% \\
		\bottomrule
	\end{tabular}%
\end{table}%

\subsection{Impact of Appointment Flexibility}
\label{sec:flexibility}

For the base case, the appointment flexibility was assumed to be $\Delta = 60$ minutes.
Deviating from a previously agreed appointment must be negotiated with the customer, but if there are significant benefits in terms of efficiency, this may be worth the effort.
To determine the impact of appointment flexibility, Table~\ref{tab:time_flexibility} presents results for the base case (small network, $\Delta=60$), in which the value of $\Delta$ is varied.
The model is given four hours of CPU time for each setting.
The columns are similar to the previous table, and show the number of automated orders, the relative and absolute savings, and the additional savings compared to the base case.
Note that the appointment flexibility does not affect the amount of orders that may benefit from automation, which is 437 for all four experiments.

Table~\ref{tab:time_flexibility} reveals that, if the appointment flexibility is already limited to one hour, limiting it further to 30 minutes to increase the service level is relatively inexpensive: the cost savings would only go down by 0.8\%.
Increasing the flexibility by 30 minutes, on the other hand, goes a long way.
Using $\Delta=90$ instead of $\Delta=60$ results in 5\% additional savings.
This indicates that the impact of appointment flexibility can be substantial.
Also note that the additional benefit is almost half that of the additional benefit for extending the network (+11\%).

\begin{table}[!t]
	\centering
	\footnotesize
	\caption{Overview of Cost Savings under Different Values of $\Delta$.}
	\label{tab:time_flexibility}%
	\begin{tabular}{crrrrr}
		\toprule
		&       &       &       &       & Additional savings \\
		\multicolumn{1}{l}{Network} & $\Delta$ & Autom. Orders & Rel. savings & Cost savings & comp. to base case \\
		\midrule
		\multirow{4}[0]{*}{Small} & 30    & 437   & 29\%  & \$ 110,887 & -0.8\% \\
		& 60    & 437   & 29\%  & \$ 111,802 & 0.0\% \\
		& 90    & 437   & 30\%  & \$ 117,354 & 5.0\% \\
		& 120   & 437   & 30\%  & \$ 116,865 & 4.5\% \\
		\bottomrule
	\end{tabular}%
\end{table}%

It is surprising to see that increasing $\Delta$ from 90 to 120 actually leads to a schedule that is less efficient, while more flexibility is available.
This is due to the CP model not finding the optimal solution.
Finding the best schedule is a very challenging task, and increasing the flexibility increases the search space, which makes this task even more challenging.
Further investigation is needed to determine the actual additional savings that can be realized for $\Delta=120$.
The results do suggest that no schedule could easily be found that was significantly better than the schedule for $\Delta=90$, which hints that the advantage of additional flexibility is leveling off after $\Delta=90$.

\section{Conclusion}
\label{sec:conclusion}

Autonomous freight transportation is expected to completely transform the industry, and the technology is advancing rapidly, with different players developing and testing high automation L4 trucks.
A crucial factor for the adoption of autonomous trucks is its return on investment, which is still uncertain.
This study contributes to the discussion by quantifying the benefits of the ATHN model, which is one of the most likely future scenarios.
The benefits are estimated based on a real transportation network, taking into account the detailed operations of the autonomous network.

A CP model was presented to schedule the orders on an ATHN, and the model was used to conduct a case study on a real transportation network.
It was found that solving this large-scale optimization problem with CP is computationally feasible.
Furthermore, ATHN may lead to substantial cost savings.
For some of the most challenging orders in the Southeast (orders that make a single delivery and return empty), operational cost may be reduced by 27\% to 40\%, which could save an estimated \$5.5M to \$8.4M per year on these orders.
This shows that the cost savings estimated by \cite{RolandBerger2018-ShiftingGearAutomation} (22\%-40\%) may indeed be realized.
The savings are mainly attributed to a reduction in labor cost, but the increased flexibility of autonomous trucks also plays a significant role: Even if autonomous trucks would have the same cost per mile as human-driven trucks, cost savings would still be possible.

It was also explored how different assumptions impact the ATHN.
Increasing the size of the autonomous network mainly increases the number of orders that can benefit from automation, while the average benefit per automated order remains similar.
The impact of the cost per mile for autonomous trucking was also studied.
As autonomous trucks become cheaper, it is cost-efficient to automate more orders, and existing trips become cheaper as well.
Due to the additional flexibility, it was found that the system benefit of automation often exceeds the benefit of the lower cost per mile.
Finally, it was analyzed how appointment flexibility impacts the efficiency, and it was found that allowing deviations to be even 30 minutes larger can go a long way.

This paper quantified the impact of autonomous trucking on a real transportation network, and substantial benefits were found in terms of labor costs and flexibility.
These results strengthen the business case for autonomous trucking, and major opportunities may arise in the coming years.
To seize these opportunities, transport operators will have to update their business models, and use optimization technology to operate the more complex systems.
Developing more detailed models and solution methods to support this transition is an interesting direction for future research.

\subsection*{Acknowledgements}

This research was funded through a gift from Ryder. Special thanks to the Ryder team for their invaluable support, expertise, and insights.

\clearpage

\bibliographystyle{splncs04}
\bibliography{DalmeijerVanHentenryck2021}

\end{document}